\documentclass[twocolumn]{amsproc}
\usepackage{amsthm,amsmath,amssymb}
\usepackage{abstract}

\usepackage{setspace}
\newtheorem{theorem}{Theorem}[section]

\newtheorem{corollary}{Corollary}[section]

\theoremstyle{remark}
\newtheorem{definition}{Definition}[section]

\newtheorem{prob}{Open Problem}[section]

\numberwithin{equation}{section}

\allowdisplaybreaks

\DeclareMathOperator{\RE}{Re}
\DeclareMathOperator{\IM}{Im}

\newcommand{\UCV}{\mathcal{U C  V}}
\newcommand{\UST}{\mathcal{U S  T}}

\newcommand{\SP}{\mathcal{S}_\mathcal{P}}
\newcommand{\USP}{\mathcal{U S  P}}
\newcommand{\UCSP}{\mathcal{U C  S  P}}
\newcommand{\TUCV}{\mathcal{T  U C  V}}
\newcommand{\TSP}{\mathcal{T S}_\mathcal{P}}
\newcommand{\TS}{\mathcal{T  S}^*}
\newcommand{\TC}{\mathcal{T C}}

\keywords{Starlike, convex, parabolic starlike, uniformly convex
functions. }

\subjclass{30C45}

\begin{document}

\title[Uniformly Convex and Uniformly Starlike Functions]
{\vspace*{.1cm} Uniformly Convex and Uniformly Starlike Functions}


\author[R. M. Ali]{\noindent Rosihan M. Ali }

\address{\newline School of Mathematical Sciences,
\newline Universiti Sains Malaysia,
  11800 USM, Penang, Malaysia }
\email{rosihan@cs.usm.my}
%
%
%

\author[V. Ravichandran]{\noindent V. Ravichandran }

\address{\newline Department of Mathematics,
\newline University of Delhi,
Delhi--110 007, India} \email{vravi68@gmail.com}

\twocolumn[ \maketitle
\begin{onecolabstract}
A normalized univalent function  is uniformly convex if it maps
every circular arc contained  in the open unit disk with center  in
it into a convex curve. This article  surveys    recent results on
the class of uniformly convex functions and on an analogous class of
uniformly starlike functions.
\end{onecolabstract}]

\section{Introduction}
One of the cornerstones in geometric function theory is the proof of
the coefficient conjecture of Bieberbach (1916) by Louis de Branges
\cite{debranges}  in the year  1985. The conjecture asserts that the
coefficient of a univalent function $f(z)=z+\sum_{n=2}^\infty a_n
z^n$ in the unit disk $\mathbb{D}= \{z\in\mathbb{C}: |z| <1\}$
satisfies $|a_n|\leq n$ with strict inequality unless $f$ is a
rotation of the Koebe function \[ k(z)=\frac{z}{(1-z)^2}.\] In fact,
de Branges proved the Milin conjecture (1971) on logarithmic
coefficients, which in turn implied the Robertson conjecture (1936)
on odd univalent functions, the Rogosinski conjecture (1943) on
subordinate functions, and finally the Bieberbach conjecture.
Milin's conjecture asserts that  the logarithmic coefficients
$\gamma_n$ of a univalent function $f(z)=z+\sum_{n=2}^\infty a_n
z^n$ defined by \[
\log\left(\frac{f(z)}{z}\right)=2\sum_{n=1}^\infty\gamma_n z^n\]
satisfy the inequality \[\sum_{k=1}^n
(n+1-k)\left(k|\gamma_k|^2-\frac{1}{k}\right)\leq 0,\quad
n=1,2,\dotsb.\] The logarithmic coefficients of the Koebe function
are $\gamma_n=1/n$ and trivially satisfy the Milin's conjecture. The
Robertson conjecture  asserts that  the inequality
\[1+|c_3|^2+\cdots+|c_{2n-1}|^2\leq n\]  is satisfied by every odd
univalent function  of the form $g(z)=z+c_3z^3+c_5z^5+\cdots$.
Rogosinski conjecture will be stated shortly. The proof that Milin
conjecture implies the other conjectures can be found in the books
on univalent functions, see for example, \cite{duren}.

The long quest for the proof of the conjecture  lead to many
profound contributions  in geometric function theory, particularly
the development of various tools for its resolution. These include
Loewner's parametric method, the area method  and Grunsky
inequalities, Milin's and FitzGerald's methods of exponentiating the
Grunsky inequalities, Baernstein's method of maximal functions, and
variational methods. Several subclasses of univalent functions were
also introduced from geometric considerations and investigated in an
attempt to settle the conjecture. Certain subclasses    are
described below.

%

 Let $\mathcal{A}$ be the class of all analytic
functions  in $\mathbb{D}$ and normalized by $f(0)=0=f'(0)-1$. Let
$\mathcal{S}$ be the subclass of $\mathcal{A}$ consisting of
univalent functions.  A domain $D$ is starlike with respect to a
point $a\in D$ if every line segment joining the point $a$ to any
other point in $D$ lies completely inside $D$. A domain starlike
with respect to the origin is  simply called starlike. A domain $D$
is convex if every line segment joining any two points in $D$ lies
completely inside $D$; in other words, the domain $D$ is convex if
and only if it is starlike with respect to every point in $D$. A
function $f\in \mathcal{S}$ is starlike if $f(\mathbb{D})$ is
starlike (with respect to the origin) while it is convex if
$f(\mathbb{D})$ is convex. The classes of all starlike and convex
functions are respectively denoted by $\mathcal{S}^*$ and
$\mathcal{C}$. Analytically, these classes are characterized by the
inequalities
\[f\in\mathcal{S}^*\Leftrightarrow \RE \left( \frac{zf'(z)}{f(z)} \right) > 0, \] and
\[ f\in\mathcal{C} \Leftrightarrow \RE \left(1+ \frac{zf''(z)}{f'(z)}
\right) > 0 .  \] More generally, for $ 0\leq \alpha<1$, let
$\mathcal{S}^* (\alpha)$ and $\mathcal{C}(\alpha)$ be the subclasses
of $\mathcal{S}$ consisting of respectively starlike functions of
order $\alpha$, and convex functions of order $\alpha$. These
classes are defined analytically by the  inequalities
\[f\in\mathcal{S}^*(\alpha)\Leftrightarrow \RE \left( \frac{zf'(z)}{f(z)} \right)
> \alpha, \]
  and \[ f\in\mathcal{C}(\alpha) \Leftrightarrow \RE \left(1+
\frac{zf''(z)}{f'(z)} \right) >\alpha .  \] Another generalization
of the class of starlike functions is   the class
$\mathcal{S}^*_\gamma$ of strongly starlike functions of order
$\gamma$, $0<\gamma\leq 1$, consisting  of $f\in\mathcal{S}$
satisfying the inequality
\[\left|\arg\left( \frac{zf'(z)}{f(z)} \right)\right|
<\frac{ \gamma \pi}{2}, \quad z\in \mathbb{D}  .\] Another related
class is the class of close-to-convex functions. A function $f\in
\mathcal{A}$ satisfying the condition
\[\RE\left( \frac{f'(z)}{g'(z)}\right) > \alpha, \quad  0\leq \alpha <1 , \]
for some (not necessarily normalized) convex univalent function $g$,
is called close-to-convex   of order $\alpha$. The class of all such
functions is denoted by $\mathcal{K}(\alpha)$. Close-to-convex
functions of order $0$ are simply called close-to-convex functions.
Using the fact that a function $f\in\mathcal{A}$ with
\[\RE(f'(z))>0\] is in $\mathcal{S}$,  close-to-convex functions can
be shown to be univalent. A function $f\in\mathcal{A}$ is starlike
with respect to symmetric points of order $\alpha$ if
 \[\RE\left( \frac{2zf'(z)}{f(z)-f(-z)}\right)  > \alpha, \quad  0\leq \alpha <1 . \]
These functions are also univalent and the class of all such
functions is denoted by $\mathcal{S}_s^*(\alpha)$. When $\alpha=0$,
this class is denoted by $\mathcal{S}_s^*$.  Coefficient estimates
for functions in all these classes can be obtained from the
coefficient estimates for functions with positive real part.

Starlikeness and convexity are hereditary properties in the sense
that every starlike (convex) function  maps each disk  $|z|<r<1$
onto a starlike (convex) domain. However,  Brown \cite{brown} showed
it is not always true that  $f\in\mathcal{S}^*$ maps each disk
$|z-z_0|<\rho< 1-|z_0|$ onto a domain starlike with respect to
$f(z_0)$. He did prove that the result is true for each
$f\in\mathcal{S}$ and for all sufficiently small disks in
$\mathbb{D}$. This motivates the definition of uniformly starlike
functions, though it was introduced independently of the work of
Brown \cite{brown}. For this purpose,  the notion of starlikeness
and convexity of curves is needed. Let $\gamma$ be a curve in
$\mathbb{D}$. Then the curve $\gamma$ is starlike with respect to
$w_0$ if $\arg( \gamma(t)-w_0) $ is a non-decreasing function of
$t$. The arc $\gamma$ is  convex if the argument of the tangent to
$\gamma(t)$ is a non-decreasing function of $t$.

\begin{definition}[{\cite[Definition 1, p.\ 364]{goodman-ust},
\cite[Definition 1, p.\ 87]{MR1145573}}] A function $f\in
\mathcal{S}$ is \textbf{uniformly starlike}  if $f$ maps every
circular arc $\gamma$ contained in $\mathbb{D}$ with center
$\zeta\in \mathbb{D}$ onto a starlike arc with respect to
$f(\zeta)$. The function $f\in \mathcal{S}$ is \textbf{uniformly
convex} if $f$ maps every circular arc $\gamma$ contained in
$\mathbb{D}$ with center $\zeta\in \mathbb{D}$  onto a convex arc.
The classes of  uniformly starlike functions and    uniformly convex
functions are denoted respectively  by $\UST$ and $\UCV$.
\end{definition}

This article  surveys  results on uniformly starlike and uniformly
convex functions. While there is quite a bit of literature on
uniformly convex functions,   not much is known about uniformly
starlike functions. The survey by  R\o nning \cite{MR1344982}
provides a summary of   early works on uniformly starlike and
uniformly convex functions.

\section{Uniformly Starlike Functions}

\subsection{Analytic characterization and basic properties}

The following two-variable analytic characterization of the class
$\UST$ is important for obtaining information about functions in the
class $\UST$.

\begin{theorem}\cite[Theorem 1, p.\
365]{goodman-ust}\label{ust-twovari} The function $ f$ is in $\UST $
if and only if
\begin{equation}
\RE \left( \frac{(z - \zeta)f'(z)}{f(z)-f(\zeta)} \right) \geq \ 0,
\quad z, \zeta \in  \mathbb{D} .
\end{equation}
\end{theorem}
By taking $\zeta=-z$ in the above theorem, evidently the class
$\UST\subset \mathcal{S}^*_s$ and hence $|a_n|\leq 1$ for $f\in
\UST$. A better bound $|a_n|\leq 2/n$ for $f\in \UST$, proved by
Charles Horowitz, was also reported in Goodman \cite[Theorem 4, p.\
368]{goodman-ust}. The proof involved showing  $\UST$ is a subclass
of  $\UST^*$ consisting of functions $f\in \mathcal{A}$ for which
$e^{i\alpha}f'(z)$ have positive real part  for some real number
$\alpha$.

\begin{prob}Determine the sharp coefficient estimates for
functions in the class  $\UST$ of uniformly starlike functions.
\end{prob}

Using Theorem~\ref{ust-twovari}, Goodman \cite{goodman-ust}  showed
that the function \[  F_1 (z)=\frac{z}{1-Az}\in \UST \Leftrightarrow
|A|\le \frac{1}{\sqrt{2}}.\] Similarly, if $F_2 (z)=z+Az^n$, $n>1$,
and
\[ |A|\le \frac{\sqrt 2 }{2n},\] then he  showed that $F_{2}$ is in
$\UST$. Merkes and Salamasi \cite{merkes} improved the bound to be
\[ |A|\leq \sqrt{\frac{n+1}{2n^3}}.\] For $n\neq 2$, the bound need not be
sharp. The sharp upper bound was obtained by Nezhmetdinov
\cite[Corollary 4, p.\ 47]{nezh1}. The class $\UST$ can also be seen
to be preserved under the transformations $e^{-i\alpha} f(
e^{i\alpha} z)$ and $f(tz)/t$, where $\alpha \in\mathbb{R}$ and $0<
t \leq 1$. For a given locally univalent analytic function $f\in
\mathcal{A}$, the disk automorphism is the function
$\Lambda_f:\mathbb{D}\rightarrow \mathbb{C}$ given by
\[\Lambda_f(z) := \frac{f(\varphi(z))-f(\lambda)}{(1-|\lambda|^2)f'(\lambda)},
\quad \varphi(z) = \frac{z+\lambda}{1+\overline{\lambda} z}. \] A
family $\mathcal{F}$ is linearly invariant if $\Lambda_f\in
\mathcal{F}$ whenever $f\in\mathcal{F}$. The families $\mathcal{S}$
of univalent functions  and $\mathcal{C}$ of convex functions are
linearly invariant families.  The disk automorphism of the function
$F_1$ with $A=1/2$ is not in $\UST$. This shows that the class
$\UST$ is not a linearly invariant family.

To provide another application of the above theorem, expand  the
function \[ \frac{(z - \zeta)f'(z)}{f(z)-f(\zeta)}\] in its Taylors
series in powers of $z$ and $\zeta$ respectively. Use of the
inequality $|c_n|\leq 2 \RE c_0$ for a function
$p(z)=c_0+c_1z+c_2z^2+\cdots$ with positive real part in
$\mathbb{D}$  yields the following result:

\begin{theorem} \cite[Lemma 1, p.\
365]{goodman-ust}\label{ust-positive} Let $f \in \UST$, and define
$p_0, p_1, q_0, q_1$ by \begin{alignat*}{2}
 p_0 (\zeta) & =\frac{f(\zeta )}{\zeta }, \quad & p_1(z)& =\frac{f(\zeta )(1-2a_2
\zeta )-\zeta }{\zeta ^2}, \\
 q_0 (\zeta)& =\frac{f(z)}{z{f}'(z)},\quad & q_1(z)&=\frac{f(z)-z}{z^2{f}'(z)}.
\end{alignat*} Then
\[ | p_1 (\zeta )| \le 2\RE (p_0 (\zeta )),\quad  \text{and}\quad
q_1 (z)|\le 2\RE ( q_0 (z)). \]
\end{theorem}

Theorem~\ref{ust-positive}  and the coefficient estimate $|a_n|\leq
2/n$ for $f\in \UST$ yield the growth inequality for $\UST$:
\[ \frac{r}{1+2r} \leq |f(z)|\leq -r +2\ln \frac{1}{1-r}, \quad |z|=r<1. \]
This inequality provides the lower bound for the Koebe constant for
the family $\UST$: \[ \frac{1}{3} \leq K(\UST) \leq
1-\frac{\sqrt{3}}{4}.
\] The upper bound  follows from the function $f$ given by
$f(z)=z+\sqrt{3}z^2/4$.

\begin{prob}Determine the sharp growth, distortion and rotation
estimates, as well as the Koebe constant for the class   $\UST$.
\end{prob}

Another application of Theorem~\ref{ust-twovari} follows from the
simple identity
\[
 \frac{f(z)-f(\zeta)}{(z - \zeta)f'(z)}= \int_0^1\frac{f'(tz+(1-t)\zeta)}
 {f'(z)}dt.
\]
Using this identity, Merkes and Salamasi \cite[Theorem 4, p.\
451]{merkes} showed that \[ f\in \UST\quad \text{ if } \quad \RE
\left(\frac{f'(w)}{f'(z)}\right)>0, \quad  z, w\in \mathbb{D}.\] If
$f\in \UST$, they also showed that \[
\RE\left(\frac{f'(w)}{f'(z)}\right)^{1/2}>0, \quad  z, w\in
\mathbb{D},\] and the exponent  1/2 is best possible.

\subsection{Convolution and Radius Problems}

The convolution (or Hadamard product) of two analytic functions
\[f(z)=z+\sum_{n=2}^\infty a_nz^n\quad \text{and}\quad g(z)=z+\sum_{n=2}^\infty
b_nz^n\] is the analytic function \[(f*g)(z):=z+\sum_{n=2}^\infty
a_nb_n z^n.\] The term ``convolution'' is used since
\[ (f*g)(z)= \frac{1}{2\pi i} \int_{|\zeta|=\rho}f\left( \frac{z}{\zeta}\right)g(\zeta)
\frac{  d\zeta}{\zeta},\quad |z|< \rho<1 . \] The classes of
starlike, convex and close-to-convex functions are closed under
convolution with convex functions. This was conjectured by  P\'olya
and Schoenberg \cite{polya} and proved by Rusche\-weyh\ and
Sheil-Small \cite{rusc}. Rusche\-weyh's
 monograph \cite{rusc2} gives a  comprehensive survey on
convolutions. To make use of this theory in the investigation of the
class $\UST$, Merkes and Salamasi \cite{merkes} proved the following
result.

\begin{theorem}[{\cite[Theorem~1, p.\ 450]{merkes}}]\label{merk-conv}
Let $f\in\mathcal{A}$. Then $f\in\UST$ if and only if  for all
complex numbers $\alpha$, $\beta$ with $|\alpha|<1$ and $|\beta|<1$,
\[ \RE\left( \frac{f(z)* \frac{z}{(1-\alpha z)(1-\beta z)} }
{ f(z)*\frac{z}{(1-\alpha z)^2} }\right)\geq 0,\quad z \in\mathbb{D}
.
\]
\end{theorem}

The following result of R\o nning \cite{MR1270313} is also useful in
using convolution technique to investigate $\UST$.

\begin{theorem}\cite[Lemma 3.3, p.\ 236]{MR1270313}\label{ust-onevari}
The function $ f\in \UST  $ if and only if
\begin{equation}
\RE \left( \frac{f(z)-f(xz)}{(1 - x)zf'(z)} \right) \geq \ 0, \quad
z \in \mathbb{D}, |x|=1 .
\end{equation}
\end{theorem}
%

Let $\mathcal{G}$ denote  the subset of $\mathcal{A}$ having the
property $\mathcal{P}$. If, for every $f \in \mathcal{F}$, $r^{-1}
f(rz)\in\mathcal{G}$ for $r\leq R$, and $R$ is the largest number
for which this holds, then   $R$ is the $\mathcal{G}$-radius (or the
radius of the property $\mathcal{P}$) in $\mathcal{F}$. Thus, the
radius of a property $\mathcal{P}$  in the set $\mathcal{F}$ is the
largest number $R$ such that every function in the set $\mathcal{F}$
has the property $\mathcal{P}$ in each disk
$\mathbb{D}_r=\{z\in\mathbb{D} : |z|<r \}$ for every $r<R$. For
example, a starlike function need not be convex; however, every
starlike function maps the disk $|z|< 2-\sqrt{3}$ onto a convex
domain and hence the radius of convexity of the class
$\mathcal{S}^*$ of starlike functions is $2-\sqrt{3}$.

Merkes and Salamasi \cite{merkes} (using Theorem~\ref{merk-conv})
and  R\o nning \cite{MR1353083} (using  Theorem~\ref{ust-onevari})
 independently showed that the $\UST$-radius of the class  $\mathcal{C}$ of
convex functions is $1/\sqrt{2}$. Merkes and Salamasi \cite[Theorem
5, p.\ 451]{merkes}   also obtained a  lower bound for the
$\UST$-radius for the class of pre-starlike functions. For $\alpha
\leq 1$, the class $\mathcal{R}_{\alpha}$ of \emph{prestarlike}
functions of order $\alpha$ consists of functions $f\in \mathcal{A}$
satisfying
\begin{align*}\begin{cases}
       f*\frac{z}{(1-z)^{2-2\alpha}}\in \mathcal{S}^{*}(\alpha),
       & \quad   \alpha<1  , \\
       \RE \frac{f(z)}{z}>\frac{1}{2}, &\quad  \alpha=1 .
       \end{cases}
     \end{align*}
Note that   $\mathcal{R}_0=\mathcal{C}$ and $\mathcal{R}_{1/2}=
\mathcal{S}^*(1/2)$.  The known radius  results  are recorded in the
following theorem.

\begin{theorem}
\begin{enumerate}\item[]
\item The $\UST$-radius for the class of univalent functions
$\mathcal{S}$ is $r_0\approx 0.3691$.

\item The $\UST$-radius $r_0^*$ for the class $\mathcal{S}^*$
satisfies \[0.369 < r_0^* \leq 1/\sqrt{7}.\]

\item The $\UST$-radius for the class of convex functions
$\mathcal{C}$ is $1/\sqrt{2}$.

\item The $\UST$-radius for the class of pre-starlike functions is
at least $(1+\alpha)/(1-\alpha)$  for
\[\frac{\sqrt{2}-1}{\sqrt{2}+1} \leq \alpha < 1.\]

\end{enumerate}
\end{theorem}

The exact value of the $\UST$-radius $r_0$ of  $\mathcal{S}$  is
obtained as the unique  root of $\varphi(t) =\pi/2$ in the interval
$[0,1]$ where $\varphi(t)$ is the expression in \cite[Equation
(2.1), p. 320]{MR1353083}.

\begin{prob}Determine the (exact)  $\UST$-radius $r_0^*$ of the class
$\mathcal{S}^*$ and  the exact $\UST$-radius of the class of
pre-starlike functions. Determine whether  the class $\UST$ is
closed under convolution with convex functions.
\end{prob}

For a given subset $\mathcal{V} \subset \mathcal{A}$, its dual set
$\mathcal{V}^*$ is defined by \[\mathcal{V}^* :=\left\{ g\in
\mathcal{A}: \frac{(f*g)(z)}{z}\neq 0 \text{ for all } f\in
\mathcal{V} \right\}.\] Nezhmetdinov \cite[Theorem 2, p.\ 43]{nezh1}
showed that the the dual set of the class $\UST$ is the subset of
$\mathcal{A}$ consisting of  functions $h:\mathbb{D}\rightarrow
\mathbb{C}$ given by
\[ h(z)= \frac{z\left(1-\frac{(w+i\alpha)}{1+i\alpha}z\right)}{(1-wz)(1-z)^2},
\quad \alpha \in\mathbb{R},\ |w|=1. \] He determined the uniform
estimate $|a_n(h)| \leq d n $ for the $n$-th Taylor coefficient of
$h$ in the dual set of $\UST$ with a sharp constant $d=\sqrt{M}
\approx 1.2557$, where $M\approx 1.5770$ is the maximum value of a
certain trigonometric expression. Using this, he   showed that
\[ \sum_{n=2}^\infty n|a_n| \leq \frac{1}{\sqrt{M}} \Rightarrow f\in\UST.\]
The bound $1/\sqrt{M}$ is sharp.

\begin{prob}  R\o nning \cite{MR1353083} proved that
$\UST\not\subset \mathcal{S}^*(1/2)$ and
posed the problem of determining the largest $\alpha$ such that
$\UST\subset \mathcal{S}^*(\alpha)$.  Nezhmetdinov \cite{nezh2}
showed that $\UST\not\subset \mathcal{S}^*(\alpha_0)$ for some
$\alpha_0\approx 0.1483$. Determine  the largest $\alpha$ such that
$\UST\subset\mathcal{S}^*(\alpha)$.
\end{prob}

\section{Uniformly Convex Functions}

\subsection{Analytic  characterizations  and parabolic starlike
functions} Recall that a  univalent function $f$ is in  the class
$\mathcal{UCV}$ of uniformly convex functions   if for every
circular arc $\gamma$ contained in $ \mathbb{D} $ with center
$\zeta\in \mathbb{D} $ the image arc $f(\gamma)$ is convex. From
this definition, the following theorem  is readily obtained.
\begin{theorem}[{\cite[Theorem 1, p.\ 88]{MR1145573}}]\label{ucv-two}
The function $ f$ belongs to $\UCV  $  if and only if
\begin{equation}\label{ucv-twovari}
\RE \left( 1 + (z -\zeta)\frac{f''(z)}{f'(z)} \right) \geq 0, \quad
z, \zeta \in  \mathbb{D}.
\end{equation}
\end{theorem}
Though the class $\mathcal{C}$ is a linear invariant family, the
class $\UCV$ is not. This was proved by Goodman \cite[Theorem 5, p.\
90]{MR1145573} by using the function \[ F(z)=\frac{z}{1-Az}.\] This
function $F\in \UCV$ if and only if $|A|\leq 1/3$.

From the geometric definition or from Theorem~\ref{ucv-two}, it is
evident that $\UCV\subset \mathcal{CV}$. However, by taking $\zeta=
-z$ in Theorem~\ref{ucv-two}, it is evident that  $\UCV\subset
\mathcal{C}(1/2)$. In view of this inclusion and the coefficient
estimate for functions in $\mathcal{C}(1/2)$, the Taylor
coefficients $a_n$ of $f\in \UCV$ satisfy $|a_n|\leq 1/n$. Unlike
the uniformly starlike functions,   uniformly convex functions admit
a  one-variable characterization, and this readily yields   several
important properties of functions in $\UCV$. This one-variable
characterization, obtained independently by R\o nning \cite[Theorem
1, p.\ 190]{MR1128729}, and  Ma and Minda \cite[Theorem 2, p.\
162]{MR1182182}, is the following result.

\begin{theorem}  \label{tha} Let $ f \in  \mathcal{A} $. Then
$ f \in \UCV $ if and only if
\begin{equation}\label{ucv-onevari}
\RE \left( 1 + \frac{zf''(z)}{f'(z)} \right) >
\left|\frac{zf''(z)}{f'(z)} \right|, \  z \in  \mathbb{D} .
\end{equation}
\end{theorem}

If $f\in \UCV$, then equation \eqref{ucv-onevari} follows from
\eqref{ucv-twovari} for a suitable choice of $\zeta$. For the
converse, the minimum principle for harmonic function is used to
restrict $\zeta$  and $z$ to $|\zeta|<|z|<1$.  With this
restriction, \eqref{ucv-twovari} immediately  follows from
\eqref{ucv-onevari}. To give a nice geometric interpretation of
\eqref{ucv-onevari}, let   \[\Omega_p := \{ w\in\mathbb{C}: \RE w>
|w -1| \}.\] The set $ \Omega_p$ is the interior of the parabola
\[(\IM w)^2=2\RE w -1\] and it is therefore symmetric with respect to
the real axis and has (1/2, 0) as its vertex. Then $f\in \UCV$ if
and only if \[ 1 + \frac{zf''(z)}{f'(z)} \in \Omega_p.\]

A class closely related to the class $\UCV$ is the class of
parabolic starlike functions defined below.

\begin{definition}\cite{MR1128729}
The class $\SP$ of parabolic starlike functions  consists of
functions $f \in \mathcal{A} $ satisfying
\[   \label{eqsp}
 \RE     \left(    \frac{zf'(z)}{f(z)}     \right)     >
\left| \frac{zf'(z)}{f(z)} - 1\right|, \quad  z\in \mathbb{D} . \]
In other words, the class $\SP$ consists of function $f=zF'$ where $
F \in \UCV$.
\end{definition}

Since the parabolic region $\Omega_p$ is contained in the half-plane
\[\{w:\RE w>1/2\}\] and the sector \[\{ w: |\arg w| <\pi/4\},\]  R\o
nning \cite{MR1128729}  noted that \[\SP\subset\mathcal{S}^*(1/2)
\cap \mathcal{S}^*_{1/2}.\]

The class $\mathcal{C}$ of convex functions and the class
$\mathcal{S}^*$ of starlike functions are connected by the Alexander
result that $ f \in  \mathcal{C}$ if and only if $ zf' \in
\mathcal{S}^*$. Such a question   between the classes $\UST$ and
$\UCV$ is in fact a   question of equality between $\UST$ and $\SP$.
It turns out (see \cite{MR1145573,MR1270313})  that that there is no
inclusion between them: \[ \UST\not\subset \SP\quad \text{  and
}\quad \SP\not\subset \UST.\]

\subsection{Examples}

 To give some examples of functions in $\UCV$
and $\SP$, note \cite{MR1981967} that \begin{equation}
\label{ucv-suffi} \left| \frac{zf''(z)}{f'(z)} \right| < \frac{1}{2}
\Rightarrow  f \in \UCV\end{equation} and\begin{equation} \left|
\frac{zf'(z)} {f(z)} -1 \right| <\frac{1}{2} \Rightarrow f \in
\SP.\end{equation} The proof    follows readily  from the
implication
\begin{align*}
|w|<\frac{1}{2}& \Rightarrow
|w|<\frac{1}{2}=1-\frac{1}{2}<1-|w|<\RE(1+w).\end{align*} A function
$ f(z)=z-\sum_{n=2}^\infty a_n z^n$ with $a_n\geq 0$ is called a
function with  negative coefficients. For   functions with negative
coefficients, the above condition is also necessary for a function
$f$ to be in $\UCV$ or $\SP$ (see \cite{MR1457247,MR1363841}). In
terms of the coefficients, the results can be stated as follows:
\begin{theorem}\label{ucv-negat} Let $f$ be a function of the form $
f(z)=z-\sum_{n=2}^\infty a_n z^n$ with $a_n\geq 0$. Then
\[ f\in\UCV \Leftrightarrow \sum_{n=2}^\infty n(2n-1)a_n\leq 1
\] and \[ f\in\SP \Leftrightarrow \sum_{n=2}^\infty
(2n-1)a_n\leq 1.\]
\end{theorem}
Denote the class of all   functions with negative coefficients  by
$\mathcal{T}$. Define \begin{alignat*}{2} \TUCV &:=\mathcal{T}\cap
\UCV,\quad &  \TSP &  :=\mathcal{T}\cap \SP,\\ \TS
&:=\mathcal{T}\cap \mathcal{S}^*, \quad \text{ and } &  \TC &  :=
\mathcal{T}\cap \mathcal{C}.\end{alignat*}  In terms of these
classes, the above result can be stated as \[ \TUCV=\TC(1/2)\quad
\text{  and }\quad  \TSP=\TS(1/2).\] For these and other related
results, see \cite{MR1457247,MR1363841}. Using
Theorem~\ref{ucv-negat}, it can be seen \cite{MR1128729} that \[
f(z) = z - A_{n}z^{n} \in \SP \Leftrightarrow |A_{n}| \leq
\frac{1}{2n-1},\] and \[f\in \UCV\Leftrightarrow |A_{n}| \leq
\frac{1}{n(2n-1)}.\] Goodman \cite{MR1145573} showed
\[\sum_{n=2}^\infty n(n-1)|a_n|\leq \frac{1}{3}\Rightarrow
f(z)=z+\sum_{n=2}^\infty a_n z^n\in \UCV;\] this easily follows from
Theorem~\ref{ucv-negat} since
\[ \sum_{n=2}^\infty n(2n-1)a_n \leq 3 \sum_{n=2}^\infty
n(n-1)|a_n|\leq 1 . \]

The sufficient condition in \eqref{ucv-suffi} can be extended to a
more general circular region. For this purpose, let $a > 1/2$. Then
it can be shown that the minimum distance from the point $w=a$ to
points on the parabola \[|w-1| = \RE w\] is given by
\[  R_{a} = \begin{cases}  a-\frac{1}{2}, & if \
\frac{1}{2} < a  \leq \frac{3}{2} \\
\sqrt{2a-2}, & if \ a \geq \frac{3}{2}.
\end{cases}  \] Thus \cite{MR1415180}
\[ \left|1+ \frac{zf''(z)}{f'(z)} -a \right| < R_a \Rightarrow  f
\in \UCV \] and \[ \left| \frac{zf'(z)} {f(z)} -a \right| <R_a
\Rightarrow f \in \SP.\]

\subsection{Subordination and its consequences} Let $f$ and $F$ be
analytic  functions in $\mathbb{D}$. Then $f$ is said to be
\emph{subordinate} to the function  $F$, written $f(z)\prec F(z)$,
if there exists an  analytic function $w: \mathbb{D} \rightarrow
\mathbb{D}$ satisfying $w(0)=0$  such that $f(z)=F(w(z))$.  If
$p:\mathbb{D}\rightarrow \mathbb{C}$, $p(0)=1$ and $\RE p(z)>0$,
then \[ p(z)\prec \frac{1+z}{1-z}.\] This follows since the mapping
$q(z)=(1+z)/(1-z)$ maps $\mathbb{D}$ onto the right-half plane
$\Omega_H:=\{w\in\mathbb{C}: \RE w>0\}$. In  this light, the classes
of starlike and convex functions can be expressed as follows:
\[
\mathcal{S}^*   = \left\{f\in\mathcal{A}: \frac{zf'(z)}{f(z)} \prec
\frac{1+z}{1-z} \right\}\] and \[  \mathcal{C}  =
\left\{f\in\mathcal{A}: 1+ \frac{zf''(z)}{f'(z)} \prec
\frac{1+z}{1-z} \right\} .\]

R\o nning \cite{MR1128729} and  Ma and Minda \cite{MR1182182} showed
that the function $\varphi_p:\mathbb{D}\rightarrow \mathbb{C}$
defined by
\begin{align}\label{parmap}
 \varphi_p(z)&=
1+\frac{2}{\pi^2}\left( \log \frac{1+\sqrt{z}}{1-\sqrt{z}} \right)^2
\\
& \notag
=1+\frac{8}{\pi^2}\left(z+\frac{2}{3}z^2+\frac{23}{45}z^3+\frac{44}{105}z^4+\cdots
\right) \end{align}  maps $\mathbb{C}$ onto the parabolic region
\[\Omega_p := \{ w\in\mathbb{C}: \RE w> |w -1| \}.\] Therefore the
classes $\UCV$ and $\SP$ can be expressed in the form
\[
\SP   = \left\{f\in\mathcal{A}: \frac{zf'(z)}{f(z)} \prec
\varphi_p(z) \right\} \] and\[  \UCV  = \left\{f\in\mathcal{A}: 1+
\frac{zf''(z)}{f'(z)} \prec \varphi_p(z) \right\} .\]

R\o nning \cite[Theorem 6, p.\ 195]{MR1128729} went on to show the
sharp inequality \[ |f'(z)|\leq   \exp\left(
\frac{14\zeta(3)}{\pi^2}\right)\approx 5.502\] for $f\in \UCV$,
where $\zeta(t)$ denotes the Riemann zeta function.   Ma and Minda
\cite{MR1182182} on the other hand obtained distortion (bounds for
$|f'(z)|$), growth (bounds for $|f(z)|$), covering (the radius of
the largest disk centered at origin contained in $f(\mathbb{D})$)
and rotation (the upper bound for $|\arg(f'(z))|$) estimates for
functions in $\UCV$. These results are then proved for more general
classes of functions by Ma and Minda \cite{MR1343506}. For this
purpose, let $\phi$ be an analytic function with positive real part
in $\mathbb{D}$,   normalized by the conditions $\phi(0)=1$ and
$\phi'(0)>0$, such that $\phi$ maps the unit disk $\mathbb{D}$ onto
a region starlike with respect to $1$ that is symmetric with respect
to the real axis. They introduced the following classes: \[
\mathcal{S}^*(\varphi): = \left\{f\in\mathcal{A}:
\frac{zf'(z)}{f(z)} \prec \varphi(z) \right\} \] and\[
\mathcal{C}(\varphi) = \left\{f\in\mathcal{A}: 1+
\frac{zf''(z)}{f'(z)} \prec \varphi(z) \right\} .\] These functions
are called  Ma-Minda starlike and convex functions respectively. For
special choices of $\varphi$, these classes become well-known
classes of starlike and convex functions. For example, for the
choice
\[\varphi_{A,B}(z)=\frac{1+Az}{1+Bz},\quad -1\leq B<A\leq 1,\] the
class $\mathcal{S}^*[A,B]:=\mathcal{S}^*(\varphi_{A,B})$ is the
class of Janowski starlike functions.  For the classes of Ma-Minda
starlike and convex functions,   the following theorem is obtained.

\begin{theorem}\label{growthdistor} \cite{MR1343506}
If $f\in \mathcal{C}(\varphi)$, then, for $|z|=r$,
\begin{alignat*}{3}  k_\varphi'(-r) & \leq & |f'(z)|  &  \leq   k_\varphi'(r), \\
-k_\varphi(-r) &\leq & |f(z)| & \leq k_\varphi(r),
\end{alignat*} where $k_\varphi:\mathbb{D}\rightarrow \mathbb{C}$ is
defined by\[ 1+ \frac{zk_\varphi''(z)}{k_\varphi'(z)} = \varphi(z).
\] Equality holds for some $z\neq 0$ if and only if $f$ is a
rotation of $k_\varphi $. Also either $f $ is a rotation of
$k_\varphi $ or $f(\mathbb{D})$ contains the disk $|w|\leq
-k_\varphi(-1)$, where
\[-k_\varphi(-1) =\lim_{r\rightarrow 1-} (-k_\varphi(-r)).\] Further, for
$|z_0|=r<1$,
\[ |\arg(f'(z_0)) | \leq \max_{|z|=r} |\arg k'_\varphi(z)|.\]
\end{theorem}

The proof relies on the subordination $f'(z)\prec k'_\varphi(z)$
satisfied by functions $f\in \mathcal{C}(\varphi)$.  Corresponding
results for functions in $\mathcal{S}^*(\varphi)$ were also obtained
by Ma and Minda \cite{MR1343506}. The distortion theorem for $f\in
\mathcal{S}^*(\varphi)$ requires some additional assumptions on
$\varphi$. Theorem~\ref{growthdistor} contains the corresponding
results for uniformly convex functions \cite{MR1182182} as special
cases. Extension of these (and other closely related) results  to
functions starlike with respect to symmetric points, conjugate
points, multivalent starlike functions, and meromorphic functions
were investigated in \cite{MR2129665,MR2566823,MR2663926}.

Let $h_\varphi:\mathbb{D}\rightarrow \mathbb{C}$ be defined by
\[ \frac{zh_\varphi'(z)}{h_\varphi(z)}=\varphi(z).\] Ma and Minda
\cite{MR1343506} proved that  \[
f\in\mathcal{S}^*(\varphi)\Rightarrow
\frac{f(z)}{z}\prec\frac{h_\varphi(z)}{z}.\]
 In the case when $\varphi$ is a convex
univalent function, this result is a special case of the following
  general result:

\begin{theorem}[Ruscheweyh {\cite[Theorem 1, p.\ 275]{MR0412407}}]
Let $\phi$ be a convex function defined in $\mathbb{D}$ with
$\phi(0)=1$. Define $F$ by
\begin{equation*}
 F(z)=z\exp\left( \int_0^z  \frac{\phi(x)-1}{x}dx \right).
 \end{equation*}
The function $f$ belongs to $ S^*(\phi)$ if and only if for all
$|s|\leq 1$ and $|t|\leq 1$, \begin{equation*}
  \frac{sf(tz)}{tf(sz)}
\prec \frac{sF(tz)}{tF(sz)}.
\end{equation*}
\end{theorem}

\begin{prob}Determine the  sharp bound of  $|f^{(n)}(z)|$
for $f\in \mathcal{C}(\varphi)$ and $f\in \mathcal{S}^*(\varphi)$.
For $f\in \mathcal{C}(\varphi)$,  the bounds for the cases $n=0,1$
are given by Theorem~\ref{growthdistor}. Similar bounds for $f\in
\mathcal{S}^*(\varphi)$ are also known with some restrictions on
$\varphi$.
\end{prob}

\subsection{Coefficient Problems}As noted earlier, the inclusion $\UCV\subset
\mathcal{C}(1/2)$ shows that each Taylor coefficient $a_n$ of $f\in
\UCV$ satisfies $|a_n|\leq 1/n$. These bounds can be improved. Since
the classes $\UCV$ and $\SP$ are connected by the Alexander relation
that $f\in \UCV$ if and only if $zf'\in \SP$, it   suffices to give
the coefficient estimate   for functions in $\SP$.

\begin{theorem}\cite[Theorem 5, p.\ 194]{MR1128729}
Let $f\in \SP$ and $f(z)=z+\sum_{n=2}^\infty a_n z^n$. Then
\begin{equation}\label{spcoeff}  |a_2|\leq c, \quad \text{ and }  \quad |a_n|\leq
\frac{c}{n-1}\prod_{k=3}^n
\left(1+\frac{c}{k-2}\right),\end{equation}  where $c=8/\pi^2$.
\end{theorem}
Let $p(z)=zf'(z)/f(z)=1+c_1z+c_2z^2+\cdots$, and $p(z)\prec
\varphi_p(z)$ where $\varphi_p$ is given by \eqref{parmap}.
Rogosinski's theorem states that $|c_k|\leq c$ for any function
$p(z)=1+c_1z+c_2z^2+\cdots$ subordinate to the convex univalent
function $P(z)=1+cz+\cdots$.    The coefficients of $f$ and the
coefficients of $p$ are related by \[(n-1)a_n =\sum_{k=1}^{n-1}
c_{n-k}a_k.\] This together with Rogosinksi's theorem yield the
desired coefficient bounds. Whenever $\varphi$ is a convex univalent
function, the bounds for $|a_n|$ for $f\in \mathcal{S}^*(\varphi)$
is also given by \eqref{spcoeff} where $c:=\varphi'(0)$. The
estimates given by \eqref{spcoeff} are not  sharp in general.
However, in the case, $\varphi(z)=(1+z)/(1-z)$, the inequalities in
\eqref{spcoeff} give sharp bounds for the coefficients of starlike
functions.

The sharp coefficient estimates for functions in $\UCV$ or $\SP$ is
still an open problem. However, the sharp  estimates of $|a_n|$ for
$f\in \UCV$ were obtained by Ma and Minda
\cite{MR1182182,MR1244399}. They    \cite[Theorem 5, p.\
172]{MR1182182}   also proved the sharp order of growth
$|a_n|=O(1/n^2)$ for $f\in \UCV$. The same order of growth holds for
$f\in \mathcal{C}(\varphi)$ if $\varphi$ belongs to the Hardy class
of analytic functions $\mathcal{H}^2$ (see \cite{MR1343506}). They
\cite{MR1244399} also found the sharp upper bound for the
Fekete-Szeg\H o functional $|\mu a_2^2-a_3|$ in the class $\UCV$ for
all real $\mu$. For the inverse function
\[f^{-1}(w)=w+\sum_{n=2}^\infty d_nw^n,\] they \cite{MR1244399}
obtained the sharp inequality \[|d_n|\leq\frac{ 8}{(n-1)n\pi^2},
\quad  n=2, 3,4 .\] More generally, the coefficient problem for
$f\in \mathcal{C}(\varphi)$ is also open. Estimates for the first
two coefficients as well as for the Fekete-Szeg\H o functional for
functions in $\mathcal{C}(\varphi)$ were obtained in
\cite{MR1343506}. For several related coefficient problems, see
\cite{MR2323552}.

\begin{theorem}\label{th1}
 Let $ \phi(z)=1+B_1z+B_2z^2+ \cdots$. If
$f(z)=z+a_2z^2+a_3z^3+\cdots\in \mathcal{C}(\varphi)$, then
\[
|a_3 - \mu a_2^2| \leq \begin{cases} \frac{1}{6}(B_2-(3/2)\mu
B_1^2+B_1^2) & \\
\quad  \text{ if } 3 B_1^2\mu \leq 2(B_2+B_1^2-B_1)\\
\frac{B_1}{6} \quad  \text{ if } 2(B_2+B_1^2-B_1) \leq 3 B_1^2\mu\\
 \quad  \leq 2(B_2+B_1^2+B_1)\\
\quad \frac{1}{6}(-B_2+(3/2)\mu B_1^2-B_1^2) & \\
\quad \text{ if }2(B_2+B_1^2+B_1)  \leq 3 B_1^2\mu\\
\end{cases}\]
 The result is sharp.
\end{theorem}
To see an outline of the proof,   first  express  the coefficient of
$f$ in terms of the coefficients $c_k$ for functions with positive
real part.   For $f\in \mathcal{C}(\varphi)$, let
$p:\mathbb{D}\rightarrow \mathbb{C}$ be defined by \[
p(z):=\frac{zf'(z)}{ f(z) } = 1+b_1z+b_2z^2+\cdots \]  so that $
2a_2= b_1 $ and $6a_3= b_2+ b_1^2$. Since $\phi$ is univalent and
$p(z)\prec \phi(z)$, the function
\[ p_1(z) = \frac{1+\phi^{-1}(p(z))}{1+\phi^{-1}(p(z))}=1+c_1z+c_2z^2+\cdots\]
is analytic and has positive real part in $\mathbb{D}$. Also
\begin{equation}\label{p5}
p(z)=\phi\left(\frac{p_1(z)-1}{p_1(z)+1}\right)
\end{equation}
and from this equation (\ref{p5}), it follows that
\[ b_1=\frac{1}{2}B_1c_1 \text{ and }
b_2=\frac{1}{2}B_1(c_2-\frac{1}{2}c_1^2)+\frac{1}{4}B_2c_1^2.\]
Therefore
\begin{equation}
 a_3-\mu a_2^2=\frac{B_1}{12}\left(c_2 - v c_1^2\right),\end{equation}
 where\[
  v:= \frac{1}{2B_1}\left(B_1-B_1^2-B_2 +\frac{3}{2}\mu
B_1^2\right).  \]  The theorem then follows by an application of the
corresponding coefficient results for function with positive real
part. Notice that this method is difficult to apply to get bounds
for $|a_n|$ for large $n$, as $a_n$ can only be expressed as a
non-linear function of the coefficients $c_k$.

\begin{prob}Determine the  sharp bound for the Taylor coefficients
$|a_n|$ $(n\geq 5)$ for $f\in \mathcal{C}(\varphi)$ and $f\in
\mathcal{S}^*(\varphi)$. The same problem for the  other classes
defined by subordination is still open.
\end{prob}

\subsection{Convolution} Recall that the \textbf{convolution}  of
two analytic functions    \[ f(z) = z +
\sum_{n=2}^{\infty}a_{n}z^{n},\quad \text{ and }\quad g(z) = z +
\sum_{n=2}^{\infty}b_{n}z^{n}\] is the analytic function  defined by
\[ (f*g)(z) := z +\sum_{n=2}^{\infty} a_{n}b_{n}z^{n}.\] The
convolution of two functions in $\mathcal{A}$ is again in
$\mathcal{A}$. Since the $n$th coefficient of normalized univalent
function is bounded by $n$, the convolution of the Koebe function
$k(z)=z/(1-z)^2$ with itself is not univalent. Thus, the convolution
of two univalent (or starlike) functions need not be univalent.
P\'olya  and Schoenberg \cite{polya} conjectured that the class of
convex functions $\mathcal{C}$ is preserved under convolution with
convex functions: \[ f,g \in \mathcal{C} \Rightarrow f*g\in
\mathcal{C}.\] In 1973, Ruscheweyh and Sheil-Small \cite{rusc} (see
also \cite{rusc2}) proved the Polya-Schoenberg conjecture. In fact,
they also proved that the classes of  starlike functions and
close-to-convex functions are closed under convolution with convex
functions. The proof of these facts follow  from the following
result which is also used below to show that the classes $\UCV$ and
$\SP$ are closed under convolution with starlike functions of order
1/2.

\begin{theorem}{\rm \cite[Theorem 2.4, p.\ 54]{rusc} }
\label{rusch-th11} Let $\alpha\leq 1$, $f\in \mathcal{R}_{\alpha}$
and $g\in \mathcal{S}^{*}(\alpha)$. Then, for any analytic function
$H\in\mathcal{H}(\mathbb{D})$,
\[\frac{f*Hg}{f*g}(\mathbb{D}) \subset \overline{co}(H(\mathbb{D})),\]
where $\overline{co}(H(\mathbb{D}))$ denote  the closed convex hull
of $H(U)$.
\end{theorem}

\begin{theorem}\cite[Theorem 3.6, p.\ 131]{MR1344982}
Let $\varphi$ be a convex function with $\RE \varphi(z) \geq
\alpha$, $ \alpha < 1$.  If $f\in\mathcal{R}_\alpha$ and $g\in
\mathcal{S}^* (\varphi)$, then $f*g\in \mathcal{S}^*(\varphi)$.
\end{theorem}

The proof of this theorem follows readily from
Theorem~\ref{rusch-th11} by putting $H(z)=zg'(z)/g(z)$.  In view of
the fact that $f\in \mathcal{C}(\varphi)$ if and only if $zf' \in
\mathcal{S}^*(\varphi)$, an  immediate consequence of the above
theorem is the corresponding result for $\mathcal{C}(\varphi)$: if
$f\in\mathcal{R}_\alpha$ and $g\in \mathcal{C} (\varphi)$, then
$f*g\in \mathcal{C}(\varphi)$ for any  convex function $\varphi$
with $\RE \varphi(z) \geq \alpha$. In particular, the classes $\UCV$
and $\SP$ are closed under convolution with starlike functions of
order 1/2. Similar results for several other related classes of
functions can be found in \cite{MR2737343,MR2665496,MR2119159}  or
references therein.

Goodman   remarked that the class $\UCV$ is preserved under the
transformation $ e^{-i \alpha} f(e^{i \alpha}z) $ and no other
transformation seems to be available. However, since  $\UCV$ is
closed under  convolution with starlike functions of order 1/2 and
in particular with convex functions,  the following result is
obtained.

\begin{corollary}\cite{MR1981967}
 Let
\begin{align*}
\Gamma_{1}(f(z)) & = zf'(z) , \\
  \Gamma_{2}(f(z))& = \frac{1}{2} [f(z) + zf'(z)]\\
  \Gamma_{3}(f(z)) & = \frac{k+1}{z^{k}}
\int_{0}^{z} \zeta^{k-1} f(\zeta) d\zeta, \ \RE  k > 0  \\
 \Gamma_{4}(f(z)) & = \int_{0}^{z} \frac{f(\zeta) - f(\eta
\zeta)}{\zeta -\eta\zeta} d\zeta, \ |\eta| \leq 1,\  \eta \not=
1.\end{align*}  Then $ \Gamma_{i}(f) \in  \UCV $ in $ |z| < r_{i} $
whenever $ f \in \UCV, $ where \[ r_{1} = \frac{1}{3}, \ r_{2} =
\frac{\sqrt{17} - 3}{2} \approx .56155,\ r_{3} = r_{4} = 1. \]
\end{corollary}

\subsection{Gaussian Hypergeometric functions} For complex  numbers
$a,b,c\in\mathbb{C}$ with $c\neq 0, -1, -2, \dotsc$, the Gaussian
hypergeometric function $F(a,b;c;z)$ is defined by the power series
\[ F(a,b,c;z):= \sum_{n=0}^\infty \frac{(a)_n(b)_n}{(c)_n}\frac{z^n}{n!}.\]
Here $(a)_0:=1 $ for $a\neq 0$ and if $n$ is a positive integer,
then $(a)_n:=a(a+1)(a+2)\cdots(a+n-1)$. For $\beta<1$ and
$\eta\in\mathbb{R}$, define the class $R_\eta(\beta)$ by \[ R_\eta
(\beta)= \left\{f\in \mathcal{A}\mid \text {Re}(e^{i\eta}
\bigl(f'(z)-\beta)\bigr)> 0\quad \text {for}\quad
z\in\mathbb{D}\right\}.\] For the Gaussian hypergeometric function
$F(a,b,c;z)$, Kim and Ponnusamy \cite{kimpon} found conditions which
would imply that $zF(a,b;c;z)$ belongs to $\UCV$ or $R_\eta(\beta)$.
Further they derived conditions under which $f\in R_\eta(\beta)$
implies \[zF(a,b;c;z)* f(z)\in \UCV.\]  In fact, by making use of
the Gauss summation theorem and Theorem~\ref{ucv-negat}, they
obtained the following sufficient condition for $zF(a,b;c;z)\in
\UCV$.

\begin{theorem}\cite[Theorem 1, p.\ 768]{kimpon}
Let $a,b\in \mathbb{C}-\{0\}$ and $c>|a|+|b|+2$.
 If
\[ \frac{\Gamma(c-|a|-|b|)\Gamma(c)}{\Gamma(c-|a|)\Gamma(c-|b|)} \times\]\[
\left( 1+\frac{2(|a|)_2(|b|)_2}{(c-2-|a|-|b|)_2}
+\frac{5|ab|}{c-|a|-|b|-1}\right) \leq 2,
\]then  $zF(a,b;c;z)\in \UCV$.
\end{theorem}
They  also obtained a weaker condition on the parameters so that the
function $zF(a,\overline{a};c;z)\in\UCV$. The following result
provides a mapping of $R_\eta(\beta)$ into $\UCV$.

\begin{theorem}\cite[Theorem 4, p.\ 771]{kimpon}
Let $a,b\in \mathbb{C}-\{0\}$ and $c>|a|+|b|+1$. If \[
2(1-\beta)\cos \eta\Bigg(
\frac{\Gamma(c-|a|-|b|)\Gamma(c)}{\Gamma(c-|a|)\Gamma(c-|b|)}\times\]\[
\left( 1+\frac{2|ab|}{c-|a|-|b|-1}\right)-1\Bigg) \leq 1,
\] and $f\in  R_\eta(\beta)$, then   $zF(a,b;c;z)*f(z)\in \UCV$.
\end{theorem}

An extension of these results to other related classes can be found,
 for example, in \cite{MR1784495,MR1739047}.

\subsection{Integral transform}
The classes $\UCV$ and $\SP$ are closed under several   integral
operators.

\begin{theorem}\cite[Theorem 1, p.\ 320]{MR1415180}
Let $ f_{i} \in \UCV $ and $\alpha_{i}$'s   be real numbers such
that $\alpha_{i} \geq 0$, and $ \sum_{1}^{n} \alpha_{i} \leq 1$.
Then the function \[ g(z) = \int_{0}^{z} \prod_{i=1}^{n}
[f_{i}'(\zeta)]^{\alpha_{i}} d\zeta \]   belongs to $\UCV$.
\end{theorem}

As an immediate consequence of this theorem,   the function $g$
defined by
\[ g(z) = \int_{0}^{z} \prod_{1}^{n} (1-A_{i} \zeta)^{-2
\alpha_{i}}d\zeta \]\[ ( 0 \leq \alpha_{i} < 1 ,\
\sum_{1}^{n}\alpha_{i} \leq 1, \  |A_{i}| \leq \frac{1}{3} , \
i=1,2, \ldots , n  )\]  belongs to $\UCV$. The first implication in
\eqref{ucv-suffi}  yields the following result.

\begin{theorem}\cite[Theorem 2, p.\ 320]{MR1415180}
If $ f \in \mathcal{A}$ satisfies \[ \left|\frac{zf'(z)}{f(z)}-
1\right|
 < \frac{1}{4},\] then
\[ g(z) = \int_{0}^{z} \left( \frac{f(\zeta)}{\zeta}\right)^{2}
d\zeta  \] belongs to $\UCV$. \end{theorem}


%

\subsection{$k$-Uniformly convex function}
Let  $k\geq 0$. A function $f\in \mathcal{S}$ is called
$k$-uniformly convex in $\mathbb{D}$ if the image of every circular
arc $\gamma$ contained in the unit disk  $\mathbb{D}$, with center
$\zeta$, $|\zeta|\leq k$, is convex. For any fixed $k\geq0$, the
class of all $k$-uniformly convex functions is denoted by $k-\UCV$.
The class $k-\UCV$ was introduced and investigated  by Kanas and
Wisinowska \cite{MR1690599}. As in the case of uniformly convex
functions,  the following theorem holds.

\begin{theorem}[\cite{MR1690599}]\label{kucv-two}
Let $f\in \mathcal{S}$. Then the following are equivalent:
\begin{enumerate}
\item  $ f\in k-\UCV $,
\item the inequality
\[
\RE \left( 1 + (z -\zeta)\frac{f''(z)}{f'(z)} \right) \geq 0
\]
holds for all $z \in  \mathbb{D}$ and for all $|\zeta|\leq k$,
\item the   inequality
\[
\RE \left( 1 + \frac{zf''(z)}{f'(z)} \right) > k
\left|\frac{zf''(z)}{f'(z)} \right| \]
 holds for all $z\in \mathbb{D}$.
\end{enumerate}
\end{theorem}

Interestingly, the class of $k$-uniformly convex functions unifies
the class of convex functions ($k=0$) and the class of uniformly
convex functions ($k=1$). Let \[\Omega_k = \{ w:\RE w > k |w-1|\}.\]
Then the region $\Omega_k$ is elliptic for $k>1$, parabolic for
$k=1$, and hyperbolic for $0<k<1$. The region $\Omega_0$ is the
right-half plane. Several properties of uniformly convex functions
extend  to $k-\UCV$ functions; these properties are treated  in
\cite{MR1690599,MR1693661,MR1693665,MR1784495,MR1875446,MR1847833}.

\subsection{Uniformly spirallike functions}

Let $\Gamma_w$ be the image of an arc $\Gamma_z: z=z(t)$, $(a\leq t
\leq b)$ under the function $f(z)$ and let $w_0$ be a point not on
$\Gamma_w$. Recall that the arc $\Gamma_w$ is starlike with respect
to $w_0$ if $\arg(w-w_0)$ is a nondecreasing function of $t$. This
condition is equivalent to
\[ \mbox{Im} \frac{f'(z)z'(t)}{f(z)-w_0}
\geq 0 \quad (a\leq t\leq b).\] The arc $\Gamma_w$ is
$\alpha$-spirallike with respect to $w_0$ if \[ \arg\frac{
z'(t)f'(z)}{f(z)-w_0}  \]  lies between $\alpha$ and $\alpha+\pi$
\cite{duren}. The function $f$ is uniformly $\alpha$-spirallike if
the image of every circular arc $\Gamma_z$ with center at $\zeta$
lying in $\mathbb{D}$ is $\alpha$-spirallike with respect to
$f(\zeta)$.  The class of all uniformly $\alpha$-spirallike
functions is denoted by $\USP(\alpha)$. Here is an analytic
description of $\USP(\alpha)$  analogous to the class $\UST$.

\begin{theorem}\cite{raviselva} \label{spiral-th2}
Let $|\alpha| <  \frac{\pi}{2}$. A function $f\in\mathcal{A}$
belongs to $\USP( \alpha)$ if and only if
\[\RE \left(e^{-i\alpha} \frac{(z-\zeta) f'(z)}{f(z)-f(\zeta)}
\right) \geq 0, \quad  z \not= \zeta,  \quad  z, \zeta \in
\mathbb{D}.
\]
\end{theorem}

The arc $\Gamma_w$ is  convex $\alpha$-spirallike if
\[\arg\left( \frac{z''(t)}{z'(t)} + \frac{ z'(t)f''(z)}{f'(z)}\right)\]
lies between $\alpha$ and $\alpha+\pi$. The function $f$ is a
uniformly convex $\alpha$-spiral function if the image of every
circular arc $\Gamma_z$ with center at $\zeta$ lying in $\mathbb{D}$
is convex $\alpha$-spirallike. The class of all uniformly convex
$\alpha$-spiral functions is denoted by $\UCSP(\alpha)$. An analytic
description of $\UCSP(\alpha)$ analogous to the class $\UCV$ is the
following:

\begin{theorem}\cite{raviselva}
Let  $f \in \mathcal{A}$. The the following are equivalent.

\begin{enumerate}
  \item  $f\in \UCSP(\alpha)$,
  \item $f$ satisfies the inequality
\[ \RE\left( e^{-i\alpha} \left( 1+ \frac{(z-\zeta)f''(z)}{f'(z)}\right)
\right)
\geq 0, \quad z\not= \zeta, \quad  z, \zeta\in \mathbb{D},  \]
  \item $f$ satisfies the inequality
  \[ \RE\left( e^{-i\alpha} \left( 1+ \frac{zf''(z)}{f'(z)}\right)\right) \geq
\left| \frac{zf''(z)}{f'(z)} \right|, \quad z \in \mathbb{D}.
\]
\end{enumerate}
\end{theorem}

For  $f \in\mathcal{A}$, define the function $s$ by \[ f'(z)=
(s'(z))^{e^{i\alpha}\cos\alpha}.\] Then $f \in \UCSP(\alpha)$ if and
only if $s \in \UCV$. In view of this connection with $\UCV$,
properties of functions in $\UCSP$ can be obtained from the
corresponding properties of $\UCV$. The classes of uniformly
spirallike and uniformly convex spirallike functions were introduced
by  Ravichandran  \emph{et al.}\ \cite{raviselva}, and for a
generalization of the class, see \cite{xu}.

\subsection{Radius problems} The determination of the radius of
starlikeness or convexity typically requires an estimate for the
real part of the quantities \[ Q_{ST}:=\frac{zf'(z)}{f(z)}\quad
\text{ and}\quad Q_{CV}:=1+\frac{zf''(z)}{f'(z)}.\] This method of
estimating the real part of $Q_{ST}$ or $Q_{CV}$ will not work for
the radius problems associated with uniformly convex functions,
parabolic starlike functions, strongly starlike functions and
several other subclasses of starlike/convex functions.  In these
cases, one need to know the region of values of  $Q_{ST}$ or
$Q_{CV}$. This idea was first used by R\o nning for computing the
sharp radius of parabolic starlikeness for univalent functions.

\begin{theorem}
The $\SP$-radius of the class $\mathcal{S}$ of univalent functions
is 0.33217 and the $\SP$-radius of the class $\mathcal{S}^*$ of
starlike functions is $1/3\approx 0.3333$ \cite[Corollary 3, Theorem
4, p.\ 192]{MR1128729}. The $\SP$-radius of the class $\mathcal{C}$
of convex functions is $1/\sqrt{2}\approx 0.7071$ \cite[Theorem 3.1
9b, p.\ 236]{MR1270313}.
\end{theorem}

The $\SP$-radii  for the following classes of  functions were
determined by Shanmugam and Ravichandran \cite{MR1415180}:
\begin{enumerate}
\item the class of close-to-starlike functions of order $\alpha$; these
are functions   $f\in\mathcal{A}$ satisfying the condition
$\RE(f(z)/g(z))>0$ for some function $g$ starlike  of order
$\alpha$.

\item the class of functions $f(z)=z+a_nz^n+\cdots$ satisfying the
condition $\RE(f(z)/z)>0$.

\item the class of functions $f\in \mathcal{A}$ satisfying
\[\left|\frac{f(z)}{g(z)} -1 \right|<1 \] for  some  function $g$ starlike of order
$\alpha$.

\end{enumerate}

R{\o}nning \cite[Theorem 4, p. 321]{MR1353083}   showed that the
radius of uniform convexity of the classes $\mathcal{S}$ and
$\mathcal{S}^*$ is $(4-\sqrt{13})/3\approx 0.1314$. Let $
\mathcal{S}^*_n[A,B]$ consists of functions
\[f(z)=z+a_{n+1}z^{n+1}+ a_{n+2}z^{n+2}+\cdots\] satisfying
\[\frac{zf' (z)}{f(z)}\prec \frac{1+Az}{1+Bz}.\] For the special
case $A=1-2\alpha$, $B=-1$, the class is denoted by
$\mathcal{S}^*_n(\alpha)$. Ravichandran, R\o nning and Shanmugam
\cite{MR1624955} investigated $\mathcal{S}^*_n(\beta)$-radius and
$\SP$-radius  for the class $\mathcal{S}^*_n[A,B]$. They   also
investigated the radii of convexity and uniform convexity in
$\mathcal{S}^*_n(0)$. Additionally  they studied the radius problems
for functions whose derivatives belong to the Kaplan classes
$\mathcal{K}(\alpha,\beta)$; their results, in special cases, yield
radius results for various classes of close-to-convex functions and
functions of bounded bound\-ary rotation. For $0\leq \alpha\leq
\beta$, the Kaplan classes $\mathcal{K}(\alpha, \beta)$ can be
defined as follows. A function $f'\in\mathcal{K}(\alpha, \beta)$ if
and only if there is a function
$g\in\mathcal{S}^*((2+\alpha-\beta)/2)$ and a real number
$t\in\mathbb{R}$ such that
\[ \left|\arg \left( e^{it} \frac{zf'(z)}{g(z)}\right) \right|\leq
\frac{\alpha \pi}{2}.\] For the  radius of uniform convexity of a
closely related class, see \cite{MR1371463} wherein they
investigated $\mathcal{S}^*(\beta)$-radius and $\SP$-radius of
certain integral transforms and Bloch functions. Related radius
results can also be found in \cite{MR2223875}.

\subsection{Neighborhood problems}

Given $\delta \geq 0$,  Rus\-che\-weyh \cite{rus1} defined the
$\delta$-\emph{neighbourhood} $N_\delta(f)$ of a function
$f(z)=z+\sum_{n=2}^\infty a_n z^n  \in \mathcal{A}$ to be the set
\begin{equation*}
N_\delta(f) := \Bigg\{ g : g(z) = z + \sum_{k=2}^\infty b_k z^k
\mbox{ and }
 \sum_{k=2}^\infty k|a_k - b_k| \leq \delta \Bigg\}.
 \end{equation*}
Ruscheweyh \cite{rus1} proved among other results that
\[N_{1/4}(f)\subset \mathcal{S}^*\] for $f\in \mathcal{C}$. For a more
general notion  of $T$-$\delta$-neigh\-bour\-hood of an analytic
function, see Sheil-Small and Silvia \cite{tss}. Padmanabhan
\cite{MR1281499} investigated the neighbourhood problem for the
class $\UCV$.  Since the class $\UCV$ is closed under convolution
with starlike functions of order 1/2, it follows that the function
$(f(z)+\epsilon z)/(1+\epsilon)\in \SP$ for $|\epsilon|<1/4$. Using
\[f\in\SP\Leftrightarrow \frac{1}{z}(f*h)(z)\neq 0,\quad
t\in\mathbb{R}, z\in \mathbb{D},\] where \[h(z):=
\frac{2}{1-2it-t^2}\left(\frac{z}{(1-z)^2}-\left(\frac{t^2+1}{2}+it\right)
\frac{z}{1-z}\right),\]  Padmanabhan proved that $N_\delta(f)\subset
\SP$ whenever \[\frac{f(z)+\epsilon z}{1+\epsilon}\in \SP\] for
$|\epsilon|<\delta<1$. These two assertions together show  that
\[N_{1/8}(f)\subset \SP\] for $f\in \UCV$. For some related results,
see \cite{MR1693665}.

\end{document}